\documentclass[twoside]{article} 
\usepackage{graphicx}
\date{} 
\title{On the $\nu$-zeros of the Bessel functions of purely imaginary order}
\author{\sc R. B.\ Paris \\
{\em Division of Computing and Mathematics,} \\
{\em Abertay University, Dundee DD1 1HG, UK}}
\begin{document}
\def\f#1#2{\mbox{${\textstyle \frac{#1}{#2}}$}}
\def\dfrac#1#2{\displaystyle{\frac{#1}{#2}}}
\def\boldal{\mbox{\boldmath $\alpha$}}
\newcommand{\bee}{\begin{equation}}
\newcommand{\ee}{\end{equation}}
\newcommand{\sa}{\sigma}
\newcommand{\ka}{\kappa}
\newcommand{\al}{\alpha}
\newcommand{\la}{\lambda}
\newcommand{\ga}{\gamma}
\newcommand{\eps}{\epsilon}
\newcommand{\om}{\omega}
\newcommand{\fr}{\frac{1}{2}}
\newcommand{\fs}{\f{1}{2}}
\newcommand{\g}{\Gamma}
\newcommand{\br}{\biggr}
\newcommand{\bl}{\biggl}
\newcommand{\ra}{\rightarrow}
\newcommand{\gtwid}{\raisebox{-.8ex}{\mbox{$\stackrel{\textstyle >}{\sim}$}}}
\newcommand{\ltwid}{\raisebox{-.8ex}{\mbox{$\stackrel{\textstyle <}{\sim}$}}}
\renewcommand{\topfraction}{0.9}
\renewcommand{\bottomfraction}{0.9}
\renewcommand{\textfraction}{0.05}
\newcommand{\mcol}{\multicolumn}
\date{}
\maketitle
\pagestyle{myheadings}
\markboth{\hfill \sc R. B.\ Paris  \hfill}
{\hfill \sc Zeros of the Bessel functions of imaginary order\hfill}
\begin{abstract}
The $\nu$-zeros of the Bessel functions of purely imaginary order are examined for fixed argument $x>0$. In the case of the modified Bessel function of the second kind $K_{i\nu}(x)$, it is known that it possesses a countably infinite sequence of real $\nu$-zeros described by $\nu_n\sim \pi n/\log\,n$ as $n\to\infty$. Here we apply a unified approach to determine asymptotic estimates of the $\nu$-zeros of the modified Bessel functions $L_{i\nu}(x)\equiv I_{i\nu}(x)+I_{-i\nu}(x)$ and $K_{i\nu}(x)$  and the ordinary Bessel functions $J_{i\nu}(x)\pm J_{-i\nu}(x)$.
Numerical results are presented to illustrate the accuracy of the expansions so obtained.
\vspace{0.3cm}

\noindent {\bf Mathematics subject classification (2020):} 33C10, 34E05, 41A30, 41A60 
\vspace{0.1cm}
 
\noindent {\bf Keywords:} modified Bessel functions, ordinary Bessel function, imaginary order, zeros, asymptotic expansion
\end{abstract}

\vspace{0.3cm}

\noindent $\,$\hrulefill $\,$

\vspace{0.3cm}

\begin{center}
{\bf 1.\ Introduction}
\end{center}
\setcounter{section}{1}
\setcounter{equation}{0}
\renewcommand{\theequation}{\arabic{section}.\arabic{equation}}
Bessel functions of purely imaginary order play an important role in quantum mechanics in the solution of the one-dimensional Schr\"odinger equation with exponential potentials.  The Bessel function $K_{i\nu}(x)$ stands alone among the standard Bessel functions
since it is real when the argument $x$ is positive and decays as $x\to+\infty$; the Bessel functions $J_{i\nu}(x)$, 
$Y_{i\nu}(x)$ and the modified Bessel function $I_{i\nu}(x)$ are all complex when $\nu$ and $x$ are real and nonzero.

A comprehensive discussion of the properties of the Bessel functions of imaginary order has been given by Dunster in \cite{TMD}. In particular, he introduced the real function \cite[(2.2)]{TMD}
\bee\label{e11}
L_{i\nu}(x)=\frac{\pi}{2\sinh \pi\nu} \{I_{i\nu}(x)+I_{-i\nu}(x)\}\qquad (\nu\neq 0),
\ee
which is an appropriate numerically satisfactory companion to $K_{i\nu}(x)$ when $\nu$ is real and nonzero and $x>0$.
This may be contrasted with the standard definition of $K_{i\nu}(x)$ given by
\bee\label{e12}
K_{i\nu}(x)=\frac{\pi i}{2\sinh \pi\nu}\{I_{i\nu}(x)-I_{-i\nu}(x)\}.
\ee
For the $J$-Bessel function of imaginary order, Dunster \cite[(3.4)]{TMD} introduced the real functions (when $x>0$) 
\bee\label{e13a}
F_{i\nu}(x)=\frac{1}{2\cosh \fs\pi\nu} \{J_{i\nu}(x)+J_{-i\nu}(x)\},
\ee
\bee\label{e13b}
G_{i\nu}(x)=\frac{1}{2\sinh \fs\pi\nu} \{J_{i\nu}(x)-J_{-i\nu}(x),\}
\ee
which clearly satisfy $F_{-i\nu}(x)=F_{i\nu}(x)$ and $G_{-i\nu}(x)=G_{i\nu}(x)$.
He obtained several properties satisfied by these functions, such as connection formulas, Wronskians, integral representations and the behaviour at singularities. In addition, he considered the zeros of $L_{i\nu}(x)$, $F_{i\nu}(x)$ and $G_{i\nu}(x)$ as a function of the real argument $x$.

In a recent paper, Bagirova and Khanmamedov \cite{BK} considered the $\nu$-zeros of $K_{i\nu}(x)$ when $x>0$ is fixed (that is, the zeros considered as a function of $\nu$ rather than the argument $x$). They have shown that $K_{i\nu}(x)$ has a countably infinite number of (simple) real zeros in $\nu$ when $x>0$ is fixed. We label the zeros $\nu_n$ ($n=1, 2 \ldots$) and observe that it is sufficient to consider only the case $\nu>0$ since $K_{-i\nu}(x)=K_{i\nu}(x)$. By transforming the differential equation satisfied by $K_{i\nu}(x)$ into a one-dimensional Schr\"odinger equation with an exponential potential, these authors employed the well-known quantisation rule  to deduce the leading asymptotic behaviour of the $n$th zero given by
\bee
\nu_n\sim \frac{\pi n}{\log\,n}\qquad (n\to+\infty).
\ee
This approximation was further discussed and refined in \cite{KR}. The $\nu$-zeros of the Hankel functions have been considered in \cite{JC, KRG}.
The numerical computation of the Bessel functions of purely imaginary order has been discussed in \cite{MF}.

A more detailed investigation of the asymptotic behaviour of the $n$th zero of $K_{i\nu}(x)$ for fixed $x>0$ has been given by the author in \cite{P}.
This relied on use of the well-known asymptotic expansion of $K_{i\nu}(x)$ for $\nu\to+\infty$ and yielded very accurate results from a three-term expansion for $\nu_n$. Here we consider the $\nu$-zeros of the functions $L_{i\nu}(x)$, $K_{i\nu}(x)$, $F_{i\nu}(x)$ and $G_{i\nu}(x)$  employing a simpler, unified approach. Although this approach does not produce results as accurate as those obtained in \cite{P}, it nevertheless yields very satisfactory agreement (from a three-term expansion) with the numerically computed $\nu$-zeros of these functions.

\vspace{0.6cm}

\begin{center}
{\bf 2.\ The modified Bessel functions}
\end{center}
\setcounter{section}{2}
\setcounter{equation}{0}
\renewcommand{\theequation}{\arabic{section}.\arabic{equation}}
The modified Bessel function $I_{i\nu}(x)$ is defined by
\bee\label{e20}
I_{i\nu}(x)=\frac{(\fs x)^{i\nu}}{\g(1+i\nu)} \sum_{k=0}^\infty \frac{(\fs x)^{2k}}{k! (1+i\nu)_k},
\ee
where $(a)_\nu=\g(a+\nu)/\g(a)$ is the Pochhammer symbol. 
We can obtain the expansion of $I_{i\nu}(x)$ for $\nu\to+\infty$ by making use of the well-known expansion
\[\frac{1}{\g(1+i\nu)}\sim \frac{1}{\sqrt{2\pi i\nu}}\,e^{i\nu} (i\nu)^{-i\nu} \sum_{k=0}^\infty \frac{\gamma_k}{(i\nu)^k}\qquad(\nu\to+\infty),\]
where the first few Stirling coefficients $\gamma_k$ are
\[\gamma_0=1,\quad \gamma_1=-\f{1}{12},\quad \gamma_2=\f{1}{288},\quad \gamma_3=\f{139}{51840},\quad \gamma_4=-\f{571}{2488320}\quad \gamma_5=-\f{163879}{209018880}.\]
We observe that the sum appearing in (\ref{e20}) for fixed $x$ can be written in the form
\[\sum_{k=0}^\infty \frac{(\fs x)^{2k}}{k! (1+i\nu)_k}\sim\sum_{k=0}^\infty \frac{C_k(\chi)}{(i\nu)^k},\qquad \chi:=\frac{x^2}{4}\qquad (\nu\to+\infty),\]
where the coefficients $C_k(\chi)$ for $0\leq k\leq 5$ are
\[C_0(\chi)=1,\quad C_1(\chi)=\chi,\quad C_2(\chi)=\frac{\chi}{2} (-2+\chi),\]
\bee\label{e2coeff}
C_3(\chi)=\frac{\chi}{6}(6-9\chi+\chi^2),\quad C_4(\chi)=\frac{\chi}{24}(-24+84\chi-24\chi^2+\chi^3),
\ee
\[C_5(\chi)=\frac{\chi}{120}(120-900\chi+500\chi^2-50\chi^3+\chi^4).\]

Then we obtain the large-$\nu$ expansion of $I_{i\nu}(x)$ in the form
\bee\label{e21}
I_{i\nu}(x) \sim \frac{e^{\frac{1}{2}\pi \nu}}{\sqrt{2\pi\nu}}\,e^{-i\Phi} \sum_{k=0}^\infty\frac{\gamma_k}{(i\nu)^k} \sum_{k=0}^\infty \frac{C_k(\chi)}{(i\nu)^k}
=\frac{e^{\frac{1}{2}\pi \nu}}{\sqrt{2\pi\nu}}\,e^{-i\Phi} \sum_{k=0}^\infty\frac{a_k}{(i\nu)^k},
\ee
where 
\bee\label{e23}
\Phi:=\nu \log\,\nu-\nu+\nu\log \fs x+\frac{\pi}{4}=\nu \log\,\la\nu+\frac{\pi}{4},\qquad \la:=\frac{2}{ex}.
\ee
The coefficients $a_k$ are 
\[a_0=1,\quad a_1=C_1(\chi)+\gamma_1 C_0(\chi),\quad a_2=C_2(\chi)+\gamma_1 C_1(\chi)+\gamma_2 C_0(\chi),\]
and in general,
\[a_k=C_k(\chi)+\sum_{r=1}^k \gamma_r C_{k-r}(\chi)\qquad (k\geq1).\]
The same procedure applied to $I_{-i\nu}(x)$ leads to the expansion
\bee\label{e21b}
I_{-i\nu}(x) \sim\frac{e^{\frac{1}{2}\pi \nu}}{\sqrt{2\pi\nu}}\,e^{i\Phi} \sum_{k=0}^\infty\frac{a_k}{(-i\nu)^k}\qquad(\nu\to+\infty).
\ee

The expansions (\ref{e21}) and (\ref{e21b}) can now be employed in the definitions of the modified Bessel functions $L_{i\nu}(x)$ and $K_{i\nu}(x)$ defined in (\ref{e11}) and (\ref{e12}). We then obtain 
\bee\label{e25a}
L_{i\nu}(x) \sim \frac{\pi}{\sinh \pi\nu}\,\frac{e^{\frac{1}{2}\pi\nu}}{\sqrt{2\pi\nu}}\bl\{\cos \Phi \bl(1-\frac{a_2}{\nu^2}+\frac{a_4}{\nu^4}-\cdots\br)-\sin \Phi \bl(\frac{a_1}{\nu}-\frac{a_3}{\nu^3}+\frac{a_5}{\nu^5}- \cdots\br)\br\}\ee
and
\bee\label{e25b}
K_{i\nu}(x) \sim \frac{\pi}{\sinh \pi\nu}\,\frac{e^{\frac{1}{2}\pi\nu}}{\sqrt{2\pi\nu}}\bl\{\sin \Phi \bl(1-\frac{a_2}{\nu^2}+\frac{a_4}{\nu^4}-\cdots\br)+\cos \Phi \bl(\frac{a_1}{\nu}-\frac{a_3}{\nu^3}+\frac{a_5}{\nu^5}- \cdots\br)\br\}\ee
for $\nu\to+\infty$.
\vspace{0.3cm}

\noindent{\bf 2.1\ \ Determination of the zeros}
\vspace{0.2cm}

We first consider the $\nu$-zeros of $L_{i\nu}(x)$ and, in (\ref{e25a}), accordingly set
\bee\label{e22}
\Phi=(n+\fs)\pi-\epsilon,
\ee
where $n$ is a large positive integer and $\epsilon$ is a small quantity. From (\ref{e25a}), the $\nu$-zeros are determined by
\bee\label{e26}
\sin \epsilon \bl(1-\frac{a_2}{\nu^2}+\frac{a_4}{\nu^4}-\cdots\bl)-\cos \epsilon \bl(\frac{a_1}{\nu}-\frac{a_3}{\nu^3}+\frac{a_5}{\nu^5}-\cdots\br)=0,\ee
whence
\[\tan \epsilon=\frac{\frac{a_1}{\nu}-\frac{a_3}{\nu^3}+\frac{a_5}{\nu^5}-\cdots}{1-\frac{a_2}{\nu^2}+\frac{a_4}{\nu^4}-\cdots}=\frac{a_1}{\nu}+\frac{a_1a_2-a_3}{\nu^3}+\frac{a_1a_2^2-a_2a_3-a_1a_4+a_5}{\nu^5}+\cdots\ .\]
Inversion of the tangent then yields
\bee\label{e27}
\epsilon=\frac{A_0}{\nu}+\frac{A_1}{\nu^3}+\frac{A_2}{\nu^5}+\cdots \ \ee
for $\nu\to+\infty$, where 
\[A_0=a_1,\quad A_1=a_1a_2-a_3-\frac{1}{3}a_1^3,\] 
\bee\label{e2coeffA}
A_2=a_1a_2^2+a_1^2a_3-a_1^3a_2-a_2a_3-a_1a_4+a_5+\frac{1}{5}a_1^5.
\ee
From (\ref{e23}) and (\ref{e22}), we finally obtain the equation describing the large-$n$ zeros of $L_{i\nu}(x)$ given by
\bee\label{e24}
\nu \log\,\la\nu=(n+\f{1}{4})\pi-\epsilon=m_+-\frac{A_0}{\nu}-\frac{A_1}{\nu^3}-\frac{A_2}{\nu^5}-\cdots ,
\ee
where for convenience we have put $m_+=(n+\f{1}{4})\pi$. 

To solve this equation we now expand $\nu\equiv \nu_n$ as 
\[\nu_n=\xi+\frac{b_0}{\xi}+\frac{b_1}{\xi^3}+\frac{b_2}{\xi^5}+\cdots ,\]
where the $b_k$ are constants to be determined and we suppose that $\xi$ is large as $n\to\infty$. Substitution in (\ref{e24}) then produces
\[\xi \log \la\xi+\frac{b_0(1+\log \la\xi)}{\xi}+\frac{b_1(1+\log \la\xi)+\fs b_0^2}{\xi^3}+\frac{b_2(1+\log \la\xi)+b_0b_1-\f{1}{6}b_0^3}{\xi^5}+\cdots\]
\[=m_+-\bl(\frac{A_0}{\xi}+\frac{A_1-A_0b_0}{\xi^3}+\frac{A_2-3A_1b_0+A_0(b_0^2-b_1)}{\xi^5}+\cdots\br)\,.
\]
Equating coefficients of like powers of $\xi$, we obtain
\bee\label{e210}
\xi\log \la\xi=m_+,
\ee
and
\[b_0=\frac{-A_0}{1+\log \la\xi},\quad b_1=\frac{A_0-A_1 b_0-\fs b_0^2}{1+\log \la\xi},\]
\[b_2=\frac{3A_1b_0-A_0(b_0^2-b_1)-A_2-b_0b_1+\f{1}{6}b_0^3}{1+\log \la\xi}~.\]

The solution of (\ref{e210}) for the lowest-order term $\xi$ can be expressed in terms of the Lambert $W$ function, which is the (positive) solution\footnote{In \cite[p.~111]{DLMF} this is denoted by Wp$(z)$.} of $W(z) e^{W(z)}=z$ for $z>0$.
Rearrangement of (\ref{e210}) shows that
\[\frac{m_+}{\xi}\,e^{m_+/\xi}=\la m_+,\]
whence
\bee\label{e33}
\xi=\frac{m_+}{W(\la m_+)}.
\ee
If we define $\chi:=\xi/m$, so that by (\ref{e31}) $1+\log \la\xi=(1+\chi)/\chi$ and introduce the coefficients
\bee\label{e3c}
B_0=\frac{-A_0}{1+\chi},\quad B_1=\frac{A_0 b_0-A_1-\fs b_0^2}{\chi^2 (1+\chi)},\quad 
B_2=\frac{3A_1b_0-A_0(b_0^2-b_1)-A_2-b_0b_1+\f{1}{6}b_0^3}{\chi^4 (1+\chi)},
\ee
then we finally have the result:
\newtheorem{theorem}{Theorem}
\begin{theorem}$\!\!\!.$\ \ 
The expansion for the $n$th $\nu$-zero of $L_{i\nu}(x)$ for fixed $x>0$ is
\bee\label{e32}
\nu_n\sim \frac{m_+}{W(\la m_+)}+\frac{B_0}{m_+}+\frac{B_1}{m_+^3}+\frac{B_2}{m_+^5}+\cdots\qquad (n\to\infty),
\ee 
where $m_+=(n+\f{1}{4})\pi$, $\la=2/(ex)$ and the coefficients $B_k$ are given in (\ref{e3c}).
\end{theorem}

The same procedure can be applied to determine the $\nu$-zeros of $K_{i\nu}(x)$. In (\ref{e25b}), we now set
\[\Phi=n\pi-\epsilon,\]
where again $n$ is a large positive integer and $\epsilon$ is a small quantity, to produce the same equation given in (\ref{e26}). Consequently, the value of $\epsilon$ is given by (\ref{e27}) with the same coefficients $A_k$.  Then
\[\nu\log\,\la\nu=(n-\f{1}{4})\pi-\epsilon=m_--\frac{A_0}{\nu}-\frac{A_1}{\nu^3}-\frac{A_2}{\nu^5}-\cdots,\]
where $m_-=(n-\f{1}{4})\pi$ and so we obtain the final result:
\begin{theorem}$\!\!\!.$\ \ The expansion for the $n$th $\nu$-zero of $K_{i\nu}(x)$ for fixed $x>0$ is
\bee\label{e39}
\nu_n\sim \frac{m_-}{W(\la m_-)}+\frac{B_0}{m_-}+\frac{B_1}{m_-^3}+\frac{B_2}{m_-^5}+\cdots\qquad (n\to\infty),
\ee 
where $m=(n-\f{1}{4})\pi$, $\la=2/(ex)$ and the coefficients $B_k$ are given in (\ref{e3c}).
\end{theorem}

The leading behaviour $\nu_n\sim \xi$ can be seen from the asymptotic expansion of $W(z)$ for $z\to+\infty$ \cite[(4.13.10)]{DLMF} 
\[W(z)\sim \log\,z-\log\log\,z+\frac{\log\log\,z}{\log\,z}+ \cdots,\]
from which it follows that
\[\frac{1}{W(z)}\sim\frac{1}{\log\,z}\bl\{1+\frac{\log\log\,z}{\log\,z}+\cdots\br\}.\]
Then, for the Bessel functions $L_{i\nu}(x)$ and $K_{i\nu}(x)$, we have the leading behaviour of the $\nu$-zeros given by 
\bee\label{e34}
\nu_n\sim\frac{m_\pm}{\log \la m_\pm}\qquad (n\to\infty),
\ee
respectively.
\vspace{0.6cm}

\begin{center}
{\bf 3.\ The Bessel functions $F_{i\nu}(x)$ and $G_{i\nu}(x)$}
\end{center}
\setcounter{section}{3}
\setcounter{equation}{0}
\renewcommand{\theequation}{\arabic{section}.\arabic{equation}}
The $J$-Bessel function of purely imaginary order is defined by
\[J_{i\nu}(x)=\frac{(\fs x)^{i\nu}}{\g(1+i\nu)} \sum_{k=0}^\infty \frac{(-)^k (\fs x)^{2k}}{k! (1+i\nu)_k},\]
from which we obtain in an analogous manner to that described in Section 2 the large-$\nu$ expansion
\bee\label{e31}
J_{\pm i\nu}(x)\sim \frac{e^{\frac{1}{2}\pi\nu}}{\sqrt{2\pi\nu}}\,e^{\mp i\Phi} \sum_{k=0}^\infty \frac{{\hat a}_k}{(\pm i\nu)^k}\qquad(\nu\to+\infty),
\ee
where $\Phi$ is defined in (\ref{e23}).
The coefficients ${\hat a}_k$ are given by 
\[{\hat a}_0=1,\quad{\hat a}_k=C_k(-\chi)+\sum_{r=1}^k \gamma_k C_{r-k}(-\chi)\qquad (k\geq 1),\]
with the first few $C_k(\chi)$ defined in (\ref{e2coeff}).  

If the expansion (\ref{e31}) is substituted into the definitions of $F_{i\nu}(x)$ and $G_{i\nu}(x)$ in (\ref{e13a}) and (\ref{e13b}), we obtain
\[\cosh \fs\pi\nu \,F_{i\nu}(x)\sim \frac{e^{\frac{1}{2}\pi\nu}}{\sqrt{2\pi\nu}}\bl\{\cos \Phi \bl(1-\frac{{\hat a}_2}{\nu^2}+\frac{{\hat a}_4}{\nu^4}+\cdots\br)-\sin \Phi \bl(\frac{{\hat a}_1}{\nu}-\frac{{\hat a}_3}{\nu^3}+\frac{{\hat a}_5}{\nu^5}-\cdots\br)\br\}
\]
and 
\[\sinh \fs\pi\nu\,G_{i\nu}(x)\sim -\frac{e^{\frac{1}{2}\pi\nu}}{\sqrt{2\pi\nu}}\bl\{\sin \Phi \bl(1-\frac{{\hat a}_2}{\nu^2}+\frac{{\hat a}_4}{\nu^4}+\cdots\br)+\cos \Phi \bl(\frac{{\hat a}_1}{\nu}-\frac{{\hat a}_3}{\nu^3}+\frac{{\hat a}_5}{\nu^5}-\cdots\br)\br\}
\] 
as $\nu\to+\infty$. Proceeding as in Section 2, we set $\Phi=(n+\fs)\pi-\epsilon$ for $F_{i\nu}(x)$ and $\Phi=n\pi-\epsilon$ for $G_{i\nu}(x)$. This produces 
\[\epsilon=\frac{A_0}{\nu}+\frac{A_1}{\nu^3}+\frac{A_2}{\nu^5}+\cdots\ ,\]
where the coefficients $A_k$ are specified in (\ref{e2coeffA}) with the $a_k$ replaced by ${\hat a}_k$. 

Omitting the details, we finally obtain the result:
\begin{theorem}$\!\!\!.$\ \ The expansion for the $n$th $\nu$-zero of $F_{i\nu}(x)$ for fixed $x>0$ is
\bee\label{e32a}
\nu_n\sim \frac{m_+}{W(\la m_+)}+\frac{B_0}{m_+}+\frac{B_1}{m_+^3}+\frac{B_2}{m_+^5}+\cdots\ \qquad (n\to\infty),
\ee
and that for $G_{i\nu}(x)$ is
\bee\label{e32b}
\nu_n\sim \frac{m_-}{W(\la m_-)}+\frac{B_0}{m_-}+\frac{B_1}{m_-^3}+\frac{B_2}{m_-^5}+\cdots\ \qquad(n\to\infty),
\ee
where $m_\pm=(n\pm\f{1}{4})\pi$ and $\la=2/(ex)$. The coefficients $B_k$ are given in (\ref{e3c}) and the $A_k$ in
(\ref{e2coeffA}) with the $a_k$ replaced by ${\hat a}_k$. 

\end{theorem}

\vspace{0.6cm}

\begin{center}
{\bf 4.\ Numerical results}
\end{center}
\setcounter{section}{4}
\setcounter{equation}{0}
\renewcommand{\theequation}{\arabic{section}.\arabic{equation}}
In this section we present results to illustrate the accuracy of the expansions in Theorems 1--3. In numerical calculations it is found more accurate to use the expression for $\xi$ in terms of the Lambert function in (\ref{e33}), rather than the asymptotic estimate, since the scale in this latter series is $\log \la m_\pm$ and so requires an extremely large value of $n$ to attain reasonable accuracy.

We present numerical results in Table 1 showing the zeros of $L_{i\nu}(x)$ and $K_{i\nu}(x)$ computed using the FindRoot command in {\it Mathematica} compared with the asymptotic values determined from the expansions (\ref{e32})
and (\ref{e39}) with coefficients $B_k$, $k\leq 2$, where $\xi$ is evaluated from (\ref{e33}). It is seen that there is satisfactory agreement with the computed zeros, even for $n=1$. However, it should be remarked that the values for the $\nu$-zeros of $K_{i\nu}(x)$ obtained in \cite{P} using the standard asymptotic expansion of $K_{i\nu}(x)$ for large $\nu$ yielded far more accurate values from a three-term expansion. For example, when $n=10$ (and $x=1$) there was agreement to 8dp, whereas in Table 1 there is agreement to only 3dp. A similar set of results is shown in Table 2 for the Bessel functions $F_{i\nu}(x)$ and $G_{i\nu}(x)$. 
\begin{table}[t]
\caption{\footnotesize{Values of the zeros of $L_{i\nu}(x)$ and $K_{i\nu}(x)$ and their asymptotic estimates when $x=1$.}}
\begin{center}
\begin{tabular}{|r|r|r||r|r|}
\hline
\mcol{1}{|c}{} & \mcol{2}{c||}{$L_{i\nu}(x)$} & \mcol{2}{c|}{$K_{i\nu}(x)$}\\
\mcol{1}{|c|}{$n$} & \mcol{1}{c|}{Zero\ $\nu_n$} & \mcol{1}{c||}{Asymptotic} & \mcol{1}{c|}{Zero\ $\nu_n$} & \mcol{1}{c|}{Asymptotic}\\
\hline
&&&&\\[-0.3cm]
1 & 3.790205 & 3.786398 & 2.962549 & 2.962961\\
2 & 5.225963 & 5.223461 & 4.534491 & 4.531277\\
3 & 6.505143 & 6.503534 & 5.879867 & 5.877888\\
4 & 7.691206 & 7.690065 & 7.107584 & 7.106243\\
5 & 8.812990 & 8.812124 & 8.258936 & 8.257949\\
10& 13.861303 & 13.860936 & 13.385883 & 13.385492\\
20& 22.620755 & 22.620598 & 22.207659 & 22.207497\\
50& 45.082187 & 45.082135 & 44.732940 & 44.732888\\
\hline
\end{tabular}
\end{center}
\end{table}
\begin{table}[t]
\caption{\footnotesize{Values of the zeros of $F_{i\nu}(x)$ and $G_{i\nu}(x)$ and their asymptotic estimates when $x=1$.}}
\begin{center}
\begin{tabular}{|r|r|r||r|r|}
\hline
\mcol{1}{|c}{} & \mcol{2}{c||}{$F_{i\nu}(x)$} & \mcol{2}{c|}{$G_{i\nu}(x)$}\\
\mcol{1}{|c|}{$n$} & \mcol{1}{c|}{Zero\ $\nu_n$} & \mcol{1}{c||}{Asymptotic} & \mcol{1}{c|}{Zero\ $\nu_n$} & \mcol{1}{c|}{Asymptotic}\\
\hline
&&&&\\[-0.3cm]
1 & 3.850274 & 3.844515 & 3.045668& 3.031436\\
2 & 5.265045 & 5.263499 & 4.581762& 4.578794\\
3 & 6.534299 & 6.534022 & 5.913240 & 5.912492\\
4 & 7.714536 & 7.714724 & 7.133494 & 7.133503\\
5 & 8.832476 & 8.832846 & 8.280167 & 8.280467\\
10& 13.872097 & 13.872514 & 13.397175 & 13.397602\\
20& 22.626541 & 22.626791 & 22.213581 & 22.213837\\
50& 45.084649 & 45.084747 & 44.735426 & 44.735525\\
\hline
\end{tabular}
\end{center}
\end{table}

\end{document}